\documentclass[english,11pt]{article}
\usepackage{amsopn,amsthm,amsfonts,amsmath,amssymb,amscd}
\usepackage{babel}
\usepackage{amssymb}
\usepackage{graphicx}

\newtheorem{theorem}{Theorem}
\newtheorem{lemma}[theorem]{Lemma}
\newtheorem{proposition}[theorem]{Proposition}
\newtheorem{corollary}[theorem]{Corollary}
\newtheorem{example}[theorem]{Example}

\newenvironment{Proof*}{\hspace{0.6cm} {\sc Proof.}}

\newcommand{\NN}{\mathbb{N}}
\newcommand{\FF}{\mathbb{F}}

\newcommand{\QQ}{\mathbb{Q}}

\newcommand{\DEF}[1]{\emph{#1}}

\renewcommand{\qedsymbol}{$\blacksquare$}

\newcommand{\adots}{\mathinner{\mskip1mu\raise1pt\vbox{\kern7pt\hbox{.}}\mskip2mu\raise4pt\hbox{.}\mskip2mu
\raise7pt\hbox{.}\mskip1mu}}

\newcommand{\sh}[1]{\mathop{\rm sh}(#1)}
\newcommand{\ch}[1]{\mathop{\rm char}(#1)}

\newcommand{\n}{{\mathcal N}}
\newcommand{\nn}{{\mathcal N}_2}
\newcommand{\nb}{{\mathcal N}_B}

\newcommand{\partition}[1]{{\mathcal P}(#1)}

\newcommand{\rp}[1]{{\mathcal R}(#1)}
\newcommand{\rpt}[2]{r(#1,#2)}
\newcommand{\ord}[1]{\mathop{\rm ord}(#1)}
\newcommand{\nil}[1]{\iota(#1)}
\newcommand{\ppot}{{\mathcal P}}
\newcommand{\bas}{{\mathcal P}_B}
\newcommand{\bask}[1]{{\mathcal P}_{B, #1}}
\newcommand{\nba}{({\mathcal N}_B,A)}

\begin{document}

\author{Polona Oblak \\
 \small{polona.oblak@fmf.uni-lj.si}
}

\date{\today}

\title{The upper bound for the index of nilpotency for a matrix commuting with a given nilpotent matrix}

\maketitle 
\parindent=0cm

\bigskip
\bigskip

{\bf Abstract.} \emph{
We study the set $\partition{\nb}$ of all possible Jordan canonical forms of nilpotent matrices commuting with a given 
nilpotent matrix $B$. We describe $\partition{\nb}$ in the special case when $B$ has only 
one Jordan block and discuss some consequences. In the general case, we find the maximal possible index of 
nilpotency in the set of all nilpotent matrices commuting with a given nilpotent matrix. We consider several examples.
}

\bigskip
\bigskip

{\small

{\bf Math. Subj. Class (2000):} 15A04, 15A21, 15A27

}

\vspace{5mm}

\bigskip
\bigskip

\section{Introduction}

\bigskip

We consider the following problem:
What are possible sizes of Jordan blocks for a pair of commuting nilpotent matrices? 
Or equivalently, for which pairs of nilpotent orbits of matrices (under similarity) there exists a pair of matrices, 
one from each orbit, that commute. The answer to the question could be considered as a generalization of Gerstenhaber--Hesselink theorem
on the partial order of nilpotent orbits \cite{cmc}.

The structure of the varieties of commuting pairs of matrices and of commuting pairs of nilpotent matrices is not yet well understood. 
It was proved by Motzkin and Taussky \cite{mt} (see also Guralnick \cite{guralnick}), that the variety of pairs of 
commuting matrices was irreducible. It was Guralnick \cite{guralnick} who showed that this is no longer the case for 
the variety of triples of commuting matrices (see also Guralnick and Sethuraman \cite{gs}, Holbrook and Omladi\v c \cite{ho}, 
Omladi\v c \cite{omladic}, Han \cite{han}, \v Sivic \cite{sivic}). Recently, it was proved that the variety of 
commuting pairs of nilpotent matrices was irreducible (Baranovsky \cite{baranovsky}, Basili \cite{basili}). 
Our motivation to study the problem is to contribute to better understanding of the structure of this variety
and which might also help in understanding the (ir)reducibility of the variety of triples of commuting matrices.

We are also motivated by the problems posed by Binding and Ko\v sir \cite{kosiMST} in the multiparameter spectral theory
and Gustafson \cite{gustafson} in the module theory over commutative rings. In both problems, certain 
pairs of commuting matrices appear. The matrices from the multiparameter spectral theory generate an algebra that 
is a complete intersection, and the matrices from the theory of modules are both functions of another matrix.

Here we initiate the study of the problem.
First we list all possible Jordan forms for nilpotent matrices commuting with a single Jordan block.
In the Section \ref{sec:digraphs} we recall the  correspondence between directed graphs and 
generic matrices (also known as Gansner--Saks Theorem), which is the main tool to prove our main result, i.e.,
to compute the maximal index of nilpotency of a nilpotent matrix commuting with a given nilpotent matrix with Jordan canonical form  
$\underline{\mu}$.
In the last section we discuss some further examples.

 \bigskip

\section{The one Jordan block case and consequences}\label{sec:general}

\bigskip

Let us denote by $\n=\n(n,\FF)$ the variety of all $n \times n$ nilpotent matrices over 
a field $\FF$ of characteristic 0 and write $\nn=\{(A,B) \in \n \times \n; AB=BA\}$. 

We follow the notations used in Basili \cite{basili} and write 
$\nb=\{A \in \n; \; AB=BA\}$ for some $B \in \n$. Suppose $\mu_1 \geq \mu_2 \geq \ldots \geq \mu_t > 0$ are
the orders of Jordan blocks in the Jordan canonical form for $B$. 
We call the partition $\underline{\mu}=(\mu_1,\mu_2,\ldots,\mu_t)$ 
the \DEF{shape} of the matrix $B$ and denote it by $\sh{B}$. 
We also write $\sh{B}=(m_1^{r_1},m_2^{r_2},\ldots,m_l^{r_l})$, where $m_1 >m_2>\ldots > m_l$.

\medskip

Let $\partition{n}$ denote the set of all partitions of $n \in \NN$
and for a subset ${\cal S} \subseteq \n(n,\FF)$ write 
$$\partition{{\cal S}}=\{\underline{\mu} \in \partition{n}; \; \underline{\mu}=\sh{A} \text{ for some } A\in {\cal S} \}\, .$$
Denote also 
$$\partition{\n_2}=\{(\sh{A},\sh{B}); \; (A,B) \in \n_2 \} \subseteq \partition{n} \times \partition{n}\, .$$
Note that $(\underline{\lambda},\underline{\mu}) \in \partition{\nn}$ if and only if
$(\underline{\mu},\underline{\lambda}) \in \partition{\nn}$, i.e. $\partition{\nn}$ is symmetric.

\medskip

It is easy to see that for each $t=1,2,\ldots,n$ there exists a uniquely defined partition 
$\rpt{n}{t}:=(\lambda_1,\lambda_2,\ldots,\lambda_t) \in \partition{n}$, such that $\lambda_1-\lambda_t \leq 1$.
It can be verified that 
$\rpt{n}{t}=\left(\left\lceil \frac{n}{t} \right\rceil^r,\left\lfloor \frac{n}{t} \right\rfloor^{t-r}\right)$. 
Denote $\rp{n}=\{\rpt{n}{t}; \; t=1,2,\ldots,n\}$.
For $\underline{\mu}=(\mu_1,\mu_2,\ldots,\mu_r)$ we define in the same fashion 
$\rpt{\underline{\mu}}{t}=(\rpt{\mu_1}{t},\rpt{\mu_2}{t},\ldots,\rpt{\mu_r}{t})$.

\bigskip

Take a matrix $B \in \n(n,\FF)$ with $\dim\ker B=1$ (i.e. $\sh{B}=(n)$). Then it is well known that any matrix commuting with
$B$ is a polynomial in $B$. By computing the lengths of the Jordan chains of $B$, we observe that $\sh{B^k}=\rpt{n}{k}$ for 
$k=1,2,\ldots,n$.	
Thus we have the following.

\bigskip

\begin{proposition}\label{thm:onejorblock}
For a matrix $B$ with $\dim\ker B=1$ it follows that $\partition{\nb}=\rp{n}$.
\hfill$\blacksquare$
\end{proposition}

\bigskip

For a sequence $(a_1,a_2,\ldots,a_k)$, $a_i \in \NN$, we write
$\ord{a_1,a_2,\ldots,a_k}=(a_{\pi(1)},a_{\pi(2)},\ldots,a_{\pi(k)})$,  where  
$a_{\pi(1)}\geq a_{\pi(2)}\geq \ldots \geq a_{\pi(k)}$ and $\pi$ is a permutation of $\{1,2,\ldots,k\}$.

For $\underline{\lambda}=(\lambda_1,\lambda_2,\ldots,\lambda_t)\in \partition{n}$ we write
\begin{align*}
 \rp{\underline{\lambda}}
              =&\{\ord{\rpt{\lambda_1}{s_1},\rpt{\lambda_2}{s_2},\ldots,\rpt{\lambda_t}{s_t}}; \; s_i=1,2,\ldots,i \; \}\, .           
\end{align*}
So, for example, $(3,3,2) \in \rp{(5,3)}$, but $(4,3,1) \notin \rp{(5,3)}$.

\bigskip

\begin{proposition} \label{thm:expjor2}
 For all $\underline{\mu} \in \partition{n}$ it follows that
 $\{(\underline{\mu}, \underline{\lambda}); \underline{\lambda} \in \rp{\underline{\mu}}\} \subseteq \partition{\nn}$.
\end{proposition}

\medskip

\begin{proof}
 Take an arbitrary $B \in \n$ with $\sh{B}=\underline{\mu}=(\mu_1,\mu_2,\ldots,\mu_t)$ 
 and pick $\underline{\lambda} \in \rp{\underline{\mu}}$.
 We want to find a matrix $A \in \nb$ such that $\sh{A}=\underline{\lambda}$.
 
 By definition, $\underline{\lambda}$ is of the form $\underline{\lambda}=\ord{\rpt{\mu_1}{k_1},\rpt{\mu_2}{k_2},
 \ldots,\rpt{\mu_t}{k_t}}$.
 Let $J_{\mu_i}$ be a $\mu_i \times \mu_i$ Jordan block 
 corresponding to eigenvalue 0. Since
 $\sh{J_{\mu_i}^{k_i}}=\rpt{\mu_i}{k_i} $, it follows that the matrix
 $A=J_{\mu_1}^{k_1} \oplus J_{\mu_2}^{k_2} \oplus \ldots \oplus J_{\mu_t}^{k_t}$
 has $\sh{A}=\underline{\lambda}$ and clearly $A \in \nb$. 
\end{proof}

\bigskip

For a matrix $A$ denote by $\nil{A}$ its index of nilpotency, i.e.
$\nil{A}=\min \{i; \; A^i=0\}$. For ${\cal S} \subseteq \n(n,\FF)$ write
$\nil{{\cal S}}=\max \{i; \; i=\nil{A} \text{ for some } A \in {\cal S}\}$.

\medskip

The following Corollaries can be easily obtained from Propositions \ref{thm:onejorblock} and \ref{thm:expjor2}.

\medskip

\begin{corollary} \label{thm:oneblock}
 If $\underline{\lambda} \in \partition{n}$, then $\underline{\lambda} \in \rp{n}$ if and only if 
 $((n),\underline{\lambda}) \in \partition{\n_2}$.
\hfill$\blacksquare$  
\end{corollary}

\medskip

\begin{corollary} \label{thm:rb1}
 For $B \in \n$ it follows that $\sh{B}\in \rp{n}$ if and only if $\nil{\nb}=n$. \hfill \qedsymbol
\end{corollary}

\medskip

\begin{example}
Note that Proposition \ref{thm:expjor2} gives us a rather large subset of partitions in $\nb$.

Take $\sh{B}=(4,3,2^2,1)$. First we directly observe that 
$(4,3,2^2,1)$, $(4,3,2,1^3)$, $(4,3,1^5)$, $(4,2^3,1^2)$, $(4,2^2,1^4)$, $(4,2,1^6)$, $(4,1^8)$,
$(3,2^4,1)$, $(3,2^3,1^3)$, $(3,2^2,1^5)$, $(3,2,1^7)$, $(3,1^9)$, $(2^5,1^2)$, $(2^4,1^4)$,
$(2^3,1^6)$, $(2^2,1^8)$, $(2,1^{10})$, $(1^{12})$ are all in $\partition{\nb}$.

Next, we see that $(4,3,2^2,1)$ is included in $\rp{\underline{\mu}}$ if
$\underline{\mu}$ is $(7,5)$, $(7,4,1)$, $(7,3,2)$, $(7,2^2,1)$, $(5,4,3)$, $(5,4,2,1)$, $(4^2,3,1)$,
$(4,3^2,2)$ or $(4,3,2^2,1)$. Thus, by symmetry, also all these partitions are in $\partition{\nb}$.

Corollary \ref{thm:oneblock} shows that $\partition{\nb} \ne \partition{12}$.
The natural question is whether the list above is the entire $\partition{\nb}$. It can be verified 
with {\tt Mathematica} that the answer is negative and that there are more partitions in $\partition{\nb}$. 
However, it would take a lot of time to compute $\partition{\nb}$ without computer.
From results in Section \ref{sec:mainthm} it follows that there exists a partition 
$(\mu_1,\mu_2,\ldots,\mu_k) \in \partition{\nb}$ such that $\mu_1=9$.
\hfill$\square$
\end{example}

\bigskip

\section{Directed graphs and the Gansner-Saks theory} \label{sec:digraphs}

\bigskip

%
%

A \DEF{digraph} is a directed graph (i.e. a graph each of whose edges are directed). 
A \DEF{path of length} $k$ in a digraph $\Gamma$ is a sequence of vertices $a_1,a_2,\ldots,a_k$ such that $(a_i,a_{i+1})$ 
is an edge in $\Gamma$ for $i=1,2,\ldots,k-1$. We do allow a path to consist of a single vertex.
For a path $\ppot$ we denote by $|\ppot|$ its length. We call $a_1$ (resp. $a_k$) the \DEF{initial} (resp. \DEF{final}) vertex
of $\ppot$.
A path $a_1,\ldots,a_k$ with $a_1=a_k$ is called a \DEF{cycle}. A digraph without cycles is called an \DEF{acyclic digraph}.
 
Let $\Gamma$ be a finite acyclic digraph on $n$ vertices, with the vertices labeled from 1 to $n$.
A \DEF{$k$--path} in $\Gamma$ is a subset of the vertices that can be partitioned into $k$ or fewer disjoint paths. Let 
$\hat{d}_k=\hat{d}_k(\Gamma)$ be the largest cardinality of a $k$--path in $\Gamma$ and define $\hat{d}_0=0$ and 
$\Delta_k=\Delta_k(\Gamma)=\hat{d}_k-\hat{d}_{k-1}$. Since $\hat{d}_k \leq \hat{d}_{k+1}$, all the $\Delta_k \, 's$
are nonnegative, so we have the infinite sequence $\Delta=\Delta(\Gamma)=(\Delta_1,\Delta_2,\ldots)$ of nonnegative integers.

\bigskip

We denote by $M(i;j)$ the entry in the $i$-th row and the $j$-th column of the matrix $M$. We say that matrices $M$ and $N$ 
of the same size \DEF{have the same pattern} if $M(i;j)=0$ if and only if $N(i;j)=0$ for $i,j=1,2,\ldots,n$.

Let $\overline{\FF}$ be a field that contains $\FF$ and has at least $n^2$ algebraically independent
transcendentals over rational numbers $\QQ$. 
A matrix $M \in \n(n,\overline{\FF})$ is called \DEF{generic} if all its nonzero entries are algebraically independent
transcendentals over $\QQ$.

\smallskip

For a nilpotent matrix $A \in \n(n,\FF)$ and a generic $M \in \n(n,\overline{\FF})$ of the same pattern it follows that 
$\nil{A} \leq \nil{M}$.

\bigskip

Given a finite acyclic digraph $\Gamma$ with $n$ vertices, we assign to it a generic
$n \times n$ matrix $M_\Gamma=[m_{ij}]$ such that $m_{ij}=0$ if $(i,j)$ is not an edge in $\Gamma$ 
and such that the rest of the entries of $M_\Gamma$ are nonzero complex numbers which are independent
transcendentals. 
Since $\Gamma$ is acyclic, $M_\Gamma$ is nilpotent.
Conversely, given a nilpotent generic matrix, it corresponds, reversing the above assignment,
to an acyclic digraph on $n$ vertices.

\begin{example}\label{ex:gs}
Consider a generic nilpotent matrix   
$M=\left[\begin{array}{cccccccccccc}
         0&0&a&0&b&0\\
         0&0&0&c&0&d\\
         0&0&0&0&0&0\\
         0&0&0&0&0&0\\
         0&0&0&e&0&f\\
         0&0&0&0&0&0
     \end{array}\right]\, .$
Its digraph $\Gamma_M$ is then equal to the digraph on the Figure \ref{fig1}. \pagebreak

       \begin{figure}[htb]
       \begin{center}
        \includegraphics[height=2.8cm, width=2.8cm]{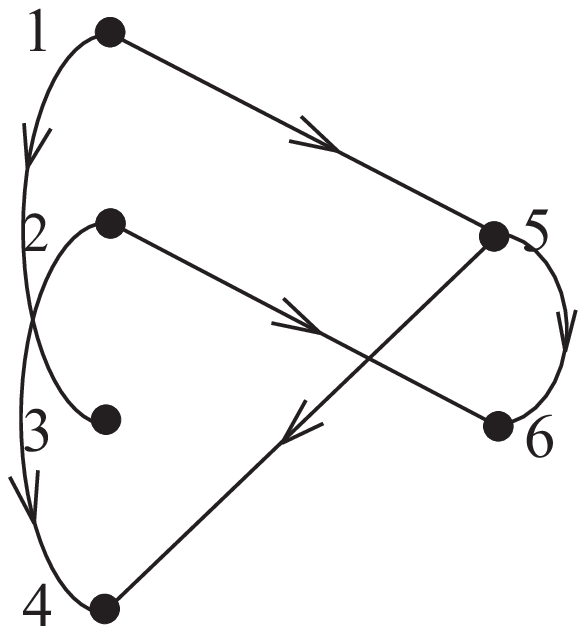}
       \end{center}
       \caption{}\label{fig1}
       \end{figure}

Conversely, every digraph $\Gamma$ as on the Figure \ref{fig1}, corresponds to a generic matrix with the
same pattern as matrix $M$. 
\hfill$\square$
\end{example}
 
 \bigskip

The following theorem was independently proved by Gansner \cite{gansner} and Saks \cite{saks}. 

\bigskip
 
\begin{theorem}\label{thm:gansnersaks}
 Let $M$ be a generic nilpotent matrix. Then 
 \begin{equation}
  \Delta(\Gamma_M)=\sh{M}\, . \tag*{$\blacksquare$}
 \end{equation}
\end{theorem}

\bigskip

\begin{corollary}\label{eq:longpath}
 The length of the longest path in an acyclic digraph $\Gamma$ is equal to $\nil{M_\Gamma}$.
\hfill$\blacksquare$
\end{corollary}
 \bigskip

\bigskip

\begin{example}
  Take matrix $M$ as in example \ref{ex:gs} and digraph $\Gamma$ as on the Figure \ref{fig1}.  
  We easily see that $\hat{d}_1=3$ (a path of length 3 is 1,5,4), $\hat{d}_2=5$ (a path of length 3 is 1,5,6 and a path of length 2 
  is 2,4) and $\hat{d}_3=6$.  By Theorem \ref{thm:gansnersaks} it follows that $\sh{M}=(3,2,1)$.
  \hfill$\square$
\end{example}

%
%

\bigskip

 \section{The $\nba$--digraph and its paths} \label{sec:preliminaries}

\bigskip

Let us fix a nilpotent matrix $B \in \n(n,\FF)$ with 
$\sh{B}=\underline{\mu}=(\mu_1,\mu_2,\ldots,\mu_t)=(m_1^{r_1},m_2^{r_2},\ldots,m_l^{r_l})$.
By convention we have $m_1 > m_2 > \ldots > m_l > 0$.

\bigskip

In this section we introduce some special digraphs, corresponding to elements in $\nb$.

\bigskip

For a pair of matrices $(A,B) \in \nn$ denote $\sh{A,B}=(\sh{A},\sh{B})$.
Then there exists $P \in GL_n(\FF)$ such that $B=P J_B P^{-1}$, where $J_B$ denotes the Jordan
canonical form of matrix $B$.
Thus $\sh{PAP^{-1},J_B}=(\sh{A},\sh{B})$ and $(PAP^{-1},J_B)\in \nn$.
Therefore we can assume that $B$ is already in its upper triangular Jordan canonical form.

\bigskip

Write $A=[A_{ij}]$ where $A_{ij} \in {\mathcal M}_{\mu_i \times \mu_j}$. It is well known 
(see e.~g.~\cite[p. 297]{GLR}) that if $AB=BA$, then $A_{ij}$ are all upper triangular Toeplitz matrices, i.e.
 for $1 \leq j \leq i \leq t$ we have
 \begin{equation}\label{eq:toeplitz}
 A_{ij}=
 \left[ \begin{matrix}
 0 & \ldots & 0 & a_{ij}^0 & a_{ij}^1 & \ldots & a_{ij}^{\mu_i-1}\\
 \vdots &  & \ddots & 0 & a_{ij}^0 & \ddots &  \vdots \\
 \vdots &  &  & \ddots & 0 & \ddots &  a_{ij}^{1} \\
 0 & \ldots & \ldots  & \ldots  & \ldots  & 0 & a_{ij}^{0}
 \end{matrix}
 \right]
 \;
 \text{ and }
 \;
 A_{ji}=
 \left[ \begin{matrix}
 a_{ji}^0 & a_{ji}^1 & \ldots & a_{ji}^{\mu_i-1}\\
 0 & a_{ji}^0 & \ddots &  \vdots \\
 \vdots &  \ddots & \ddots &  a_{ji}^{1} \\
 \vdots &   & 0 &  a_{ji}^{0} \\
 \vdots &   &  &  0 \\
 \vdots &   &  &  \vdots \\
 0 & \ldots & \ldots   & 0
 \end{matrix}
 \right]
 \, .
 \end{equation}

If $\mu_{i}=\mu_j$ then we omit the rows or columns of zeros in $A_{ij}$ or $A_{ji}$ above.

\bigskip

We introduce some further notation following Basili \cite{basili}. For a matrix $B$ with $\sh{B}=(\mu_1,\mu_2,\ldots,\mu_t)$ we
denote by $r_B$ and $k_i$, $i=1,2,\ldots,r_B$, the numbers such that 
$k_1=1$, $\mu_{k_i}-\mu_{k_{i+1}} \geq 2$, $\mu_{k_i}-\mu_{k_{i+1}-1} \leq 1$ for $i=1,2,\ldots,r_B-1$, 
$\mu_{k_{r_B}}-\mu_{t} \leq 1$.  Note that $k_i \in \{1,2,\ldots,t\}$.
For example, if $\sh{A}=(5^3,3^2,1^3)$ and $\sh{B}=(5^3,4,3^2,2^4,1^3)$, then $r_A=r_B=3$.
Furthermore, $\sh{C}\in \rp{n}$ if and only if $r_C=1$.

Set $q_0=0$, $q_l=t$ and for $\alpha=1,2,\ldots,l$ let 
$q_\alpha \in \{1,2,\ldots,t\}$, be such that $\mu_i=\mu_{i+1}$ if $q_{\alpha-1}+1 \leq i < q_\alpha$ and
$\mu_{q_\alpha} \ne \mu_{q_\alpha+1}$. For a block matrix with blocks as in
\eqref{eq:toeplitz} we define $\overline{A}_{\alpha \alpha}=[a_{ij}^0]$ where $q_{\alpha-1}+1 \leq i,j \leq q_\alpha$.

\medskip

\begin{lemma}\label{thm:basiliaa}  \cite[Proposition 2.3]{basili} 
For an $n \times n$ matrix $A$ such that $AB=BA$ it follows that $A \in \nb$ if and only if
$\overline{A}_{\alpha\alpha}$ are nilpotent for $\alpha=1,2,\ldots,l$.

Moreover, if $A \in \nb$ there exists a Jordan basis for $B$ such that $\overline{A}_{\alpha\alpha}$
are all strictly upper triangular.
\hfill{$\blacksquare$}
\end{lemma}

\bigskip

Therefore, when we study what are possible shapes of pairs from the set $\n_2$ (or in particular what 
are the indices of nilpotency in $\nb$), we may consider only matrices $A \in \nb$
such that we have $a_{ij}^0=0$ for all $1 \leq j \leq i \leq t$ whenever $\mu_i=\mu_j$.
From now, we assume that the latter relations hold for $A$.

\bigskip

Let $\Gamma$ be a digraph corresponding to a generic matrix $M \in \n(n,\overline{\FF})$ such that $M$ has 
the same pattern as a matrix $A \in \nb$. We call $\Gamma$ an \DEF{$\nba$-digraph}.

\bigskip

\begin{example} 
 Consider a nilpotent matrix $B$ with $\sh{B}=(4,2)$ and a matrix
 $A=\left[\begin{array}{cccccccccccc}
         0&0&a&0&b&0\\
         0&0&0&a&0&b\\
         0&0&0&0&0&0\\
         0&0&0&0&0&0\\
         0&0&0&c&0&d\\
         0&0&0&0&0&0
     \end{array}\right] \in \nb \, .$ 
 Then the digraph on the Figure \ref{fig1} is an $\nba$-digraph.
 \hfill$\square$
\end{example}

\bigskip

Let $\Gamma$ be an $(\nb,A)$-digraph that corresponds to $A \in \nb$. 
Denote its vertices by $(x,y)$, where $x=1,2,\ldots,t$ and $y=1,2,\ldots,\mu_x$.
Write blocks of $A$ as in  \eqref{eq:toeplitz}. If $a_{ij}^k \ne 0$ for some $1 \leq i \leq j \leq t$, 
then $\Gamma$ contains edges $\left((i,h),(j,h+k)\right)$ for $h=1,2,\ldots,\mu_j-k$.
If $a_{ij}^k \ne 0$ for some $1 \leq j \leq i \leq t$, 
then $\Gamma$ contains edges $\left((i,h),(j,h+k+\mu_j-\mu_i+1)\right)$ for all $h=1,2,\ldots,\mu_i-k$.

\bigskip

We say that edges $((i_1,j_1),(i_2,j_2))$ and $((i_3,j_3),(i_4,j_4))$ are \DEF{parallel} if 
$i_1=i_3$, $i_2=i_4$ and $j_1+j_4=j_2+j_3$. We call two paths $\ppot_1$ and $\ppot_2$ \DEF{parallel} if
they consist of pairwise parallel edges.

\bigskip

\begin{example} \label{ex:43221a}
Consider the case $\sh{B}=(4,3,2^2,1) \in \partition{12}$. Each $\nba$-digraph $\Gamma$ consists of 12 vertices:

        \begin{figure}[htb]
         \begin{center}
          \includegraphics[height=3.6cm, width=6.9cm]{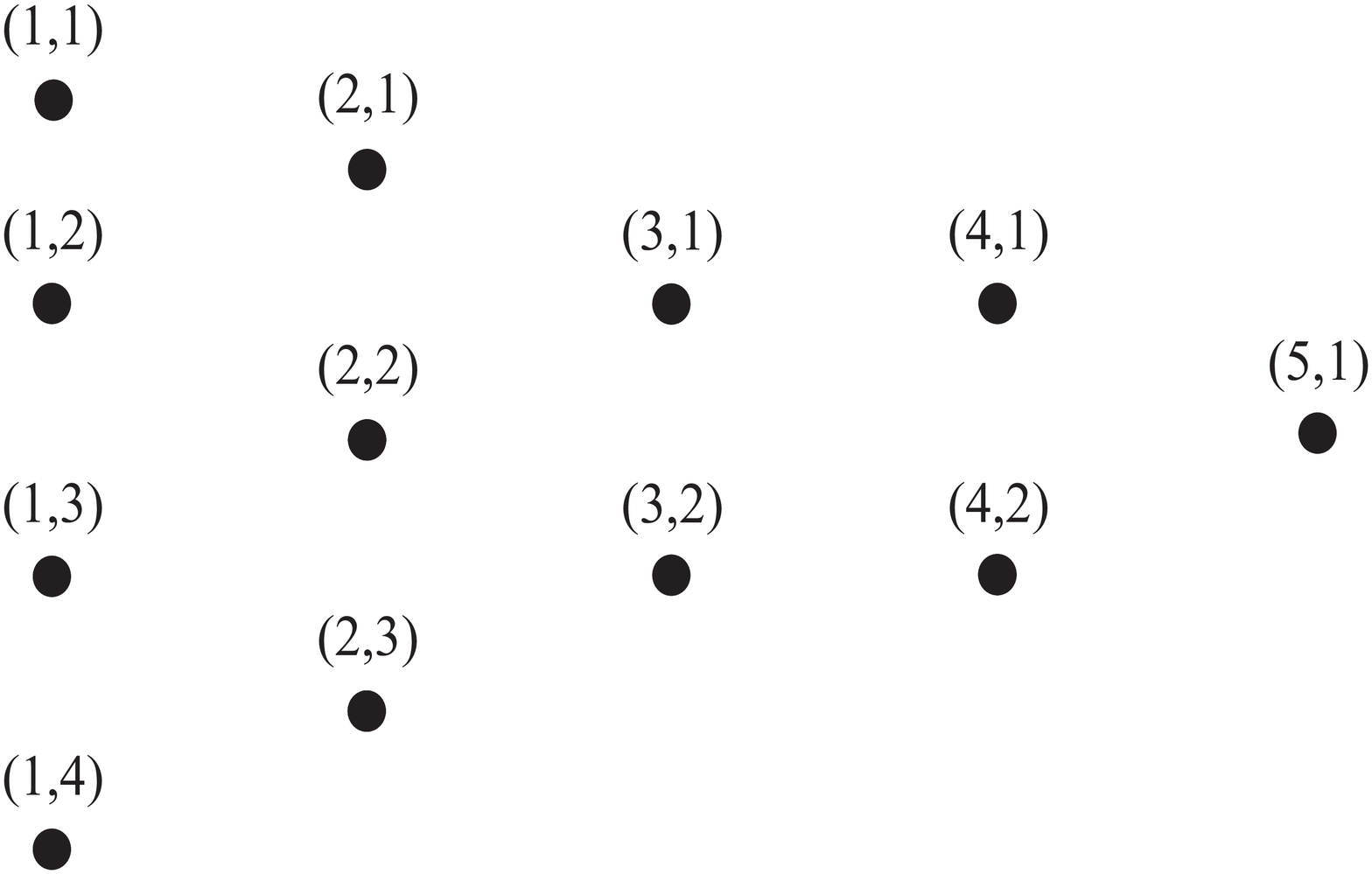}
         \end{center}
         \end{figure}

Denote by $M$ the generic nilpotent matrix, which corresponds to $\Gamma$ and has the same pattern as $A \in \nb$.
Suppose, for example, that there exists an edge $((4,1),(2,2))$ in $\Gamma$. By the correspondence 
between nilpotent matrices and acyclic digraphs it follows that $M_{10,6} \ne 0$. Since $A$ and 
$M$ have the same pattern, $a_{4,2}^0 \ne 0$ and thus
$M_{11,7} \ne 0$. Therefore, there exists also an edge $((4,2),(2,3))$ in $\Gamma$, which is parallel to 
an edge $((4,1),(2,2))$.

Similarly, we can add some other parallel edges to $\Gamma$ as described above. For example, a possible $\nba$-digraph $\Gamma$ is

        \begin{figure}[htb]
         \begin{center}
          \includegraphics[height=4.5cm, width=8.7cm]{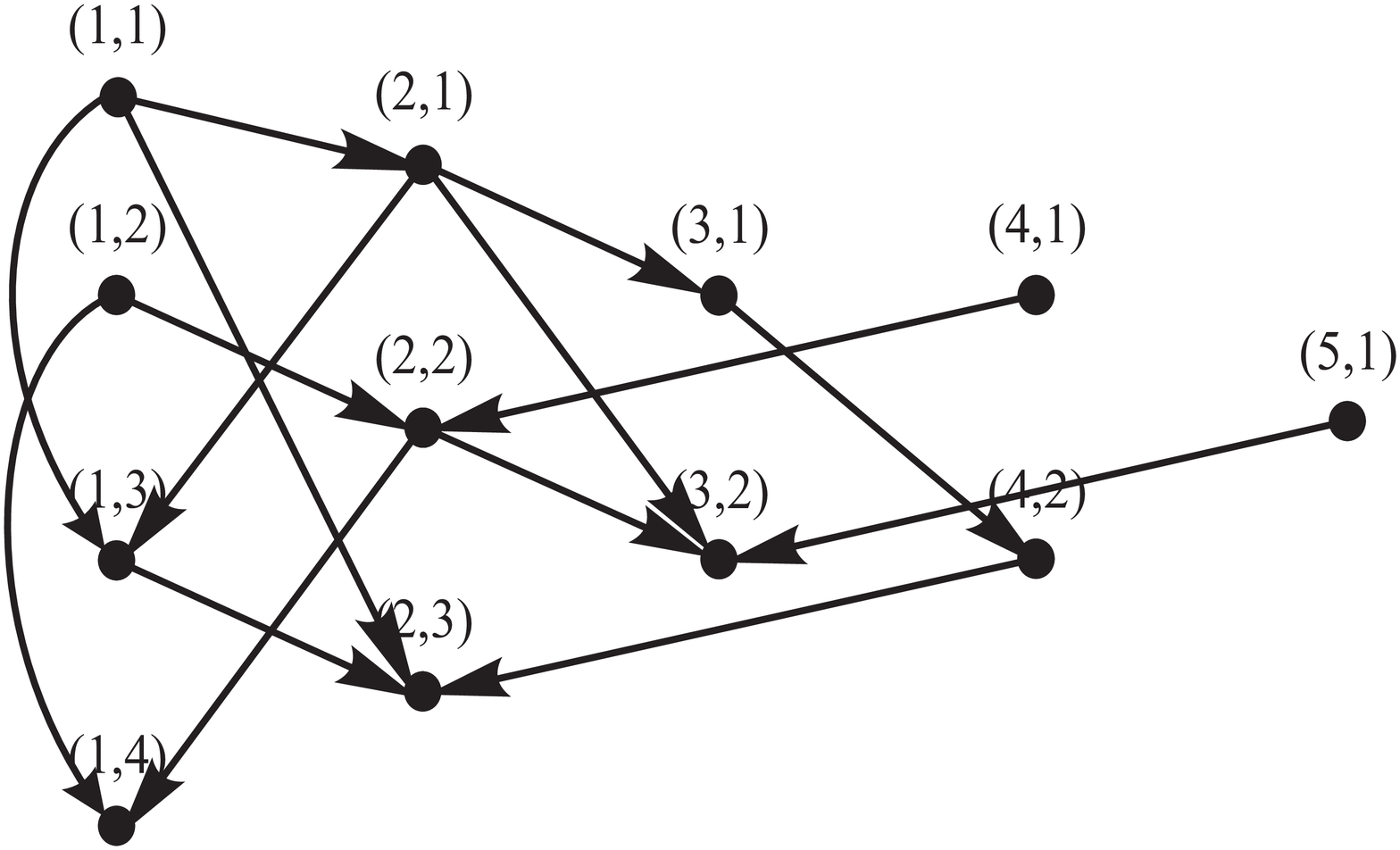}
         \end{center}
         \caption{}\label{fig3}         
         \end{figure}

and the corresponding matrix $A \in \nb$ is equal to
$$
A=\left[\begin{array}{cccc|ccc|cc|cc|c}
0 & 0 & a & 0 & b & 0 & c & 0 & 0 & 0 & 0 & 0  \\
0 & 0 & 0 & a & 0 & b & 0 & 0 & 0 & 0 & 0 & 0  \\
0 & 0 & 0 & 0 & 0 & 0 & b & 0 & 0 & 0 & 0 & 0  \\
0 & 0 & 0 & 0 & 0 & 0 & 0 & 0 & 0 & 0 & 0 & 0  \\
\hline
0 & 0 & d & 0 & 0 & 0 & 0 & e & f & 0 & 0 & 0  \\
0 & 0 & 0 & d & 0 & 0 & 0 & 0 & e & 0 & 0 & 0  \\
0 & 0 & 0 & 0 & 0 & 0 & 0 & 0 & 0 & 0 & 0 & 0  \\
\hline
0 & 0 & 0 & 0 & 0 & 0 & 0 & 0 & 0 & 0 & g & 0  \\
0 & 0 & 0 & 0 & 0 & 0 & 0 & 0 & 0 & 0 & 0 & 0  \\
\hline
0 & 0 & 0 & 0 & 0 & h & 0 & 0 & 0 & 0 & 0 & 0  \\
0 & 0 & 0 & 0 & 0 & 0 & h & 0 & 0 & 0 & 0 & 0  \\
\hline
0 & 0 & 0 & 0 & 0 & 0 & 0 & 0 & i & 0 & 0 & 0  
\end{array}\right]
$$
for some $a,b,c,\ldots,i \in \FF$.
\hfill$\square$
\end{example}

\bigskip
 
Suppose that $r_i,r_{i+1} \ne 0$ for some $i$.
Then parts $m_i^{r_i}$, $m_{i+1}^{r_{i+1}}$ in $\sh{B}$ have to be treated differently if 
$m_i-m_{i+1}=1$ or $m_{i}-m_{i+1} \geq 2$. We introduce some notation that will help us 
unify the treatments.

\bigskip

For $i=1,2,\ldots,l-1$ such that $r_i \ne 0$  define 
$$s(i)=\left\{
 \begin{array}{cl}
  r_i, & \text{if } m_{i}-m_{i+1}\geq 2 \\
  r_i+r_{i+1}, & \text{if } m_{i}-m_{i+1} \leq 1
 \end{array}
 \right.$$
and $s(l)=r_l$. 

\medskip

We set $\mu_0=\mu_1+1$. Choose any $k$, $1 \leq k \leq t$, such that $\mu_{k-1}> \mu_k > 1$ 
and let $w$ and $z$ be such that $\mu_k-\mu_w \leq 1$, $\mu_k-\mu_{w+1} \geq 2$ and 
$\mu_{k}=\mu_z > \mu_{z+1}$.
If $\mu_{k-1}> \mu_k = 1$, then set $w=z=t$.

\medskip

We denote by $V_{B,k}$ the set of vertices
\begin{align*}
 \big\{(x,y); \;  & x=k,k+1,\ldots,w, y=1,2,\ldots,\mu_{x} 
             \big\} \, .
\end{align*}                 
We call the path that contains vertices  \label{Bpath}
$$V_{1,k}=\{(x,1)\; x=1,2,\ldots,k \}, $$
$$ V_{B,k}$$ 
and 
$$ V_{3,k}=\{(x,\mu_x)\; x=1,2,\ldots,k-1 \} \cup \{(z,\mu_z)\}$$ 
the \DEF{$B_k$-path} (or \DEF{$B$-path} for short). We call $s(k)$ the \DEF{width} of $B_k$-path (or of the set $V_{B,k}$).

\pagebreak

\begin{example}\label{ex:43221c}
Consider again the case $\sh{B}=(4,3,2^2,1) \in \partition{12}$ and an $\nba$-digraph $\Gamma'$

        \begin{figure}[htb]
         \begin{center}
          \includegraphics[height=4cm, width=6cm]{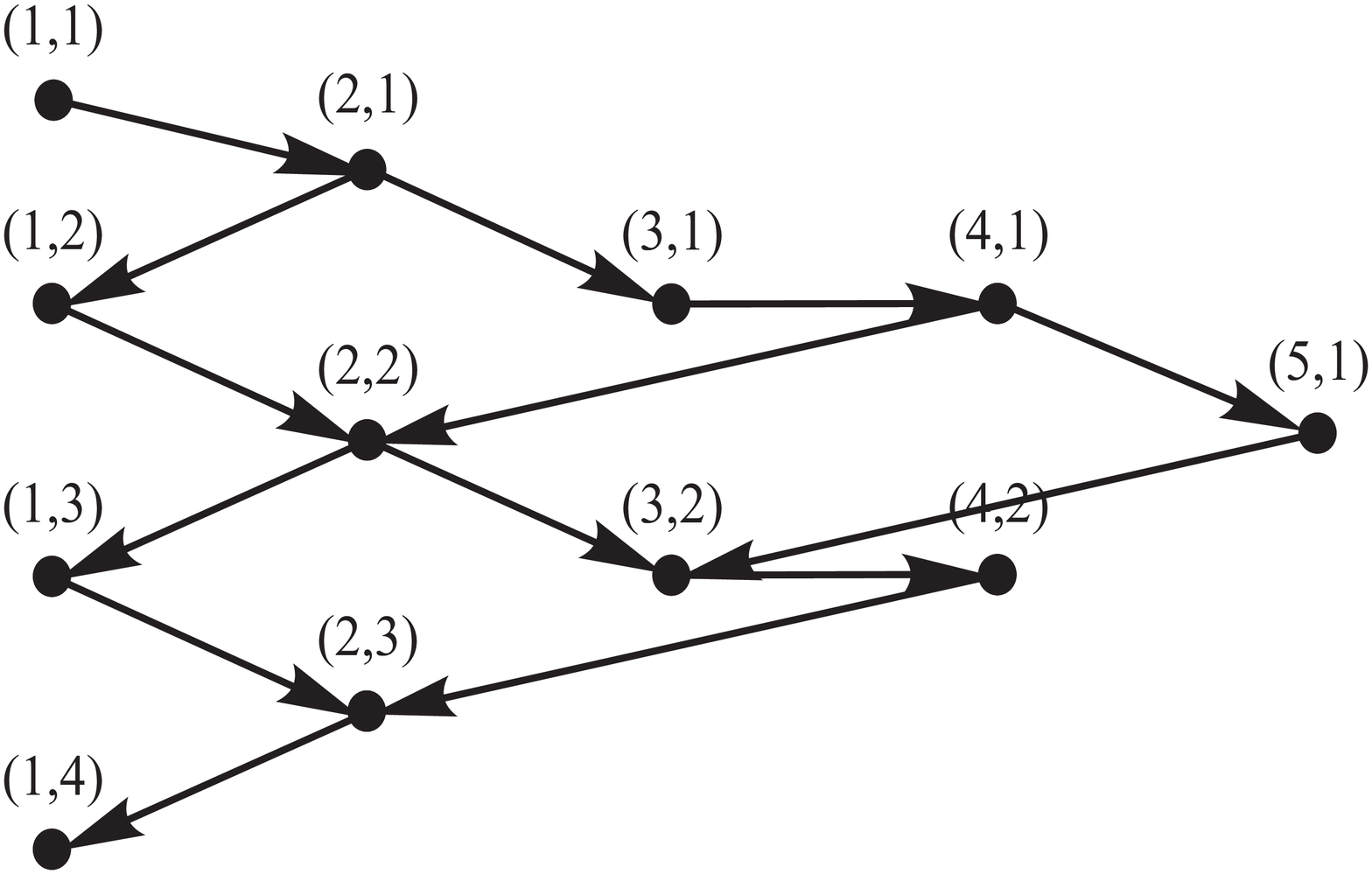}
         \end{center}
         \caption{}\label{fig4}    
         \end{figure}
         
Then for $k=1,2,3,5$, the set $V_{B,k}$ and $B_k$-path consist of the following vertices

\begin{center}
\begin{tabular}{|c||c|c|c|p{5cm}|p{6cm}|}
  \hline
 $k$ & $s(k)$ & $w$ & $z$ & $V_{B,k}$ & vertices of $B_k$-path \\
  \hline
  \hline
 $1$ & $2$ & $2$ & $1$ & $(1,1)$, $(1,2)$, $(1,3)$, $(1,4)$, $(2,1)$, $(2,2)$, $(2,3)$ & $(1,1)$, $(1,2)$, $(1,3)$, $(1,4)$, $(2,1)$, $(2,2)$, $(2,3)$ \\   
  \hline
 $2$ & $3$ & $4$ & $2$ & $(2,1)$, $(2,2)$, $(2,3)$, $(3,1)$, $(3,2)$, $(4,1)$, $(4,2)$ & $(1,1)$, $(2,1)$, $(2,2)$, $(2,3)$, $(3,1)$, $(3,2)$, $(4,1)$, $(4,2)$, $(1,4)$ \\   
  \hline
 $3$ & $3$ & $5$ & $4$ & $(3,1)$, $(3,2)$, $(4,1)$, $(4,2)$, $(5,1)$ & $(1,1)$, $(2,1)$, $(3,1)$, $(3,2)$, $(4,1)$, $(4,2)$, $(5,1)$, $(1,4)$, $(2,3)$ \\   
  \hline
 $5$ & $1$ & $5$ & $5$ & $(5,1)$ & $(1,1)$, $(2,1)$, $(3,1)$, $(4,1)$, $(5,1)$, $(1,4)$, $(2,3)$, $(3,2)$, $(4,2)$\\   
  \hline
\end{tabular}
\end{center}
\hfill$\square$
\end{example}

\bigskip

The $B_k$-path consists of three subpaths $\ppot_i$, $i=1,2,3$, such that 
$\ppot_1 \cap \ppot_2=\{(k,1)\}$ and $\ppot_2 \cap \ppot_3 = \{(z,\mu_z)\}$.
The vertices of $\ppot_1$ are from $V_{1,k}$ and its edges are 
$$((i,1),(i+1,1))$$  
for $i= 1,2, \ldots, k-1$. The $\ppot_2$ consists of vertices $V_{B,k}$ and
edges $$((i,j),(i+1,j))$$ 
for 
$i=k,k+1,\ldots,w-1$, $j=1,2,\ldots,\mu_{k}-1$ or 
$i=k,k+1,\ldots,z-1$, $j=\mu_{k}$ and 
$$((w,j),(k,j+1))$$ 
for $j=1,2,\ldots,\mu_{k}-1$.
The subpath $\ppot_3$ on vertices $V_{3,k}$ consists of edges
$$((i,\mu_{i}),(i+1,\mu_{i+1}))$$ 
if $\mu_i=\mu_{i+1}$, $1 \leq i \leq k-2$, and
$$((j,\mu_{j}),(i,\mu_{i}))$$ 
if $\mu_{k} \leq \mu_{j+1}< \mu_j < \mu_i < \mu_{i-1}$ and     
such that if $i \leq h \leq j$, then it follows that either $\mu_{h}=\mu_i \text{ or } \mu_h=\mu_j$.

\smallskip

For a $B_k$-path $\ppot$ it follows that $|\ppot|=2(k-1)+\mu_{k}+\mu_{k+1}+\ldots+\mu_{w}$.

\bigskip

\begin{example} \label{ex:43221b}
Recall the Examples \ref{ex:43221a} and \ref{ex:43221c} where $\sh{B}=(4,3,2^2,1)$. In the $\nba$-digraph $\Gamma$ from
Figure \ref{fig3} there are no $B_k$-paths. However the $\nba$-digraph $\Gamma'$ from Figure \ref{fig4}
has three $B$-paths, namely, $B_1$-path (with width 2 and length 7), $B_2$-path (with width 3 and length 9) and $B_3$-path (also of length 9
and width 3), as shown below

        \begin{figure}[htb]
         \begin{center}
          \includegraphics[height=3.2cm, width=12cm]{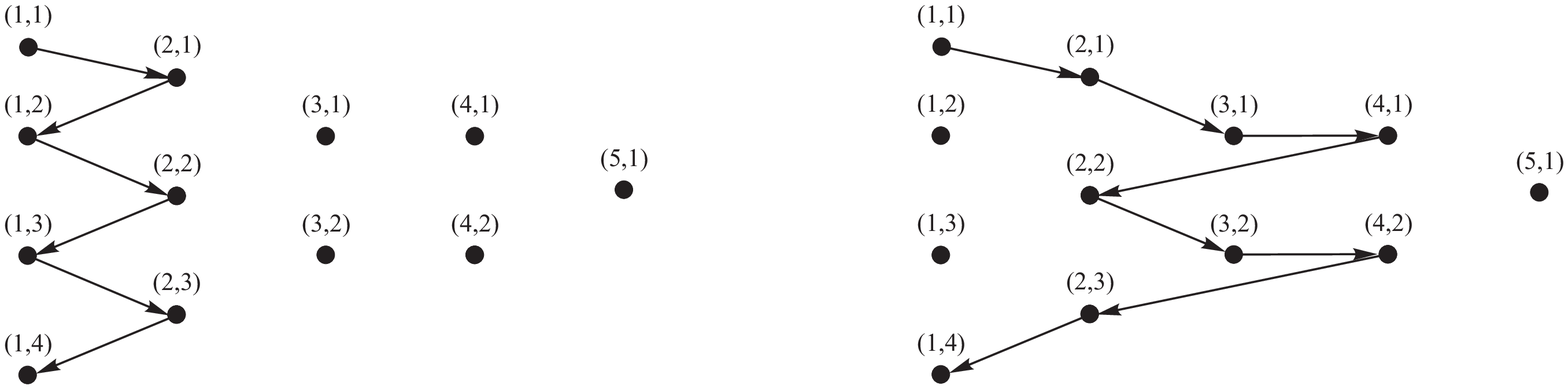}
         \end{center}
         \end{figure}

        \begin{figure}[htb]
         \begin{center}
          \includegraphics[height=3.2cm, width=4.7cm]{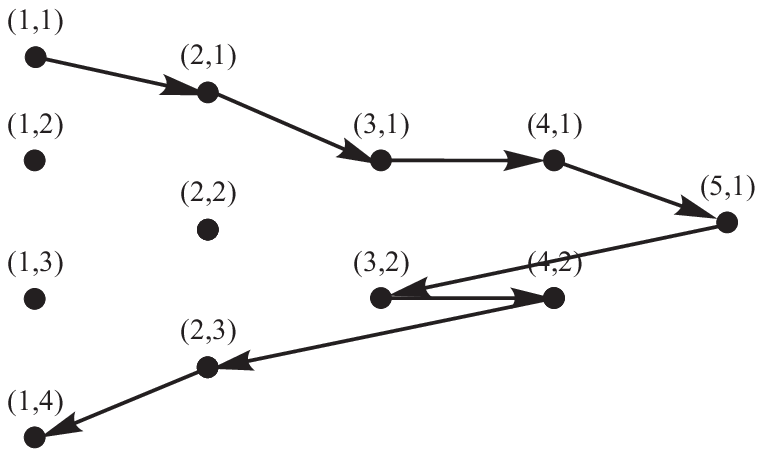}
         \end{center}
         \end{figure}

Note that in this case $B_3$-path and $B_5$-path coincide.
\hfill$\square$
\end{example}

\bigskip

\begin{example}
Let $\sh{B}=(4,3^2)$. By Corollary \ref{thm:rb1} it follows that $\nil{\nb}=10$. 
The generic $M$, such that its digraph $\Gamma_M$ is equal to

        \begin{figure}[h]
         \begin{center}
          \includegraphics[height=3.5cm, width=3.5cm]{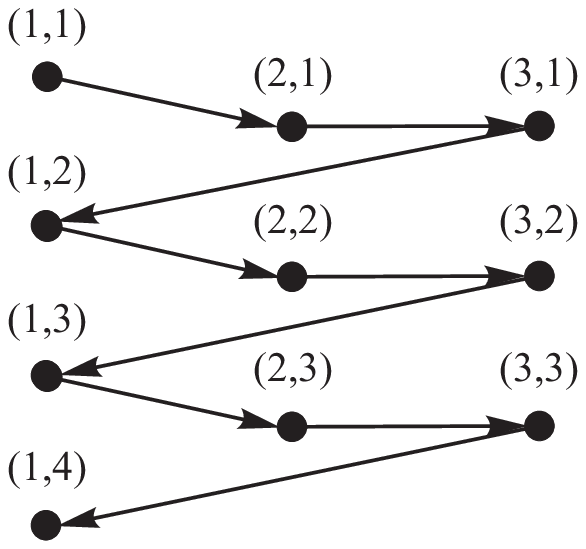}
         \end{center}
         \end{figure}

\pagebreak
and $A \in \nb$ such that $a_{12}^0$,  $a_{23}^0$, $a_{31}^0$ are its only nonzero entries,
have the same pattern. It is easy to verify that  $\nil{M}=\nil{A}=10$. 
The shown path is the only $B$-path in $\Gamma_M$.
\hfill$\square$
\end{example}

\bigskip

Recall that $\nil{A} \leq \nil{M}$ for $A \in \nb$ and a generic matrix $M \in \n(n,\overline{\FF})$ with the same pattern. 
By Corollary \ref{eq:longpath} it follows that
$\nil{\nb}$ is less than or equal to the length of the longest path in a $\nba$-digraph.
The length of the $B_k$-path is equal to 
$2 \sum_{i=1}^{k-1} r_i + r_k m_k+r_{k+1}m_{k+1}$ if
$s(k)=r_k+r_{k+1}$ and  $2 \sum_{i=1}^{k-1} r_i + r_k m_k$ otherwise.

\bigskip
  
\section{Maximal index of nilpotency in $\nb$} \label{sec:mainthm}

\bigskip

Here, we prove the following theorem:

\bigskip

\begin{theorem} \label{thm:main}
 Let $B$ be a nilpotent matrix with $\sh{B}=(\mu_1,\mu_2,\ldots,\mu_t)$. Then
  $$\nil{\nb}= \max\limits_{1 \leq i < t} \{2i+\mu_{i+1}+\mu_{i+2}+\ldots+\mu_{i+r}; \; \mu_{i+1} - \mu_{i+r}\leq 1, \mu_i \ne \mu_{i+1}\} 
  \, .$$
\end{theorem}

\bigskip

First, we show that the longest path in an $\nba$-digraph $\Gamma$ is actually equal
to the length of the longest $B$-path in $\Gamma$.

\bigskip

\begin{lemma}\label{thm:path2}
 Suppose that for $x=1,2,\ldots,t$ vertices $(x,1)$ are contained in a path $\ppot$ of an $\nba$-digraph $\Gamma$. 
 Then there exists a $B$-path $\bas$ in $\Gamma$ such that $|\ppot| \leq |\bas|$.
\end{lemma}
  
\medskip
 
\begin{proof}
 Let us write $\sh{B}=(\mu_1,\mu_2,\ldots,\mu_t)=(m_1^{r_1},m_2^{r_2},\ldots,m_l^{r_l})$, 
 where $m_i=\mu_1-i+1$ for $i=1,2,\ldots,l$ and $m_l=\mu_t$. (Note that $r_i$ can be 0 
 for some $i$.) 
 
 Let $\Gamma$ be an $\nba$-digraph and $k$ be the maximal index such that 
 $s(k)=\max\{s(i); \; 1 \leq i \leq l \}$. 
 By $\bask{k}$ we denote the $B_k$-path, i.e. the widest $B$-path in $\Gamma$. 
 
 Suppose first that $s(k)=r_k+r_{k+1}$, where $r_{k+1} \ne 0$. Then
  $$|\bask{k}|=2 \sum_{i=1}^{k-1} r_i + r_k m_k +r_{k+1}m_{k+1}=
 2 \sum_{i=1}^{k-1} r_i + (r_k+r_{k+1}) m_{k+1} + r_k \, .$$

 Examine the edge $((i,y_i),(j,y_j))$ of $\ppot$, where $\mu_i \ne \mu_j$. 
 If $i>j$  (and thus $\mu_i < \mu_j$) then  $\mu_i-y_i > \mu_j-y_j$.
 If $i<j$ (and thus $\mu_i > \mu_j$) then $y_i < y_j$. 
  
 Let $\ppot'$ be the subpath of $\ppot$ that contains all the vertices $(x,y)$ from $\ppot$, where $y > 1$.
 Let $(x_0,y_0)$ be the initial vertex of $\ppot'$. Since $((t,1),(x_0,y_0))$, $1 \leq x_0 \leq t$, is an edge of $\ppot$,
 it follows that $\mu_t-1 \geq \mu_{x_0}-y_0$.
 
 Since $\bask{k}$ is the widest $B$-path in $\Gamma$,
 it follows that
  $|\ppot'| \leq   (\mu_{x_0}-y_0) (r_{k}+r_{k+1}) + \sum_{i=1}^{l-1} r_i \leq  \sum_{i=1}^{l-1} r_i + (m_l-1) (r_{k}+r_{k+1})$. 
 Therefore
  $$|\ppot| \leq  2 \sum_{i=1}^{l-1} r_i +r_l +(m_l-1) (r_{k+1}+r_k)\, .$$
  
 Since $m_i=\mu_1-i+1$ for $i=1,2,\ldots,l$, it follows that $m_i - m_l = l-i$ for $i=1,2,\ldots,l$. 
 In particular, $m_l=m_k-l+k$. Thus
  $$|\ppot| \leq  2 \sum_{i=1}^{l-1} r_i +r_l +(m_k-l+k-1) (r_k+r_{k+1})\, .$$

 By definition of $k$ it follows that $r_i+r_{i+1} \leq r_k+r_{k+1}-1$ for all 
 $i=k+1,k+2,\ldots,l-1$.
 Thus 
 $r_{k+1}+2r_{k+2}+2r_{k+3}+\ldots+2r_{l-1}+r_l=\sum_{i=k+1}^{l-1} (r_i+r_{i+1}) \leq (l-k-1)(r_k+r_{k+1}-1)$.
 Therefore
  \begin{equation*}
   \begin{split}
   |\ppot| &\leq 2 \sum_{i=1}^{k-1} r_i + r_k + (r_k+r_{k+1}) + (l-k-1)(r_k+r_{k+1}-1)+\\
           & \qquad \qquad + (m_k-l+k-1) (r_k+r_{k+1}) =\\
     &= 2 \sum_{i=1}^{k-1} r_i + r_k + (m_{k}-1)(r_k+r_{k+1})+k+1-l \leq \\
     &\leq 2 \sum_{i=1}^{k-1} r_i + r_k + m_{k+1}(r_k+r_{k+1})\; = \;  |\bask{k}| \, .
   \end{split}
   \end{equation*}

 If $s(k)=r_k$ and $k < l$, then $r_i+r_{i+1} \leq r_k-1$ for $ i=k+1,k+2,\ldots,l-1$.
 Thus 
 $r_{k+1}+2r_{k+2}+2r_{k+3}+\ldots+2r_{l-1}+r_l=\sum_{i=k+1}^{l-1} (r_i+r_{i+1}) \leq (l-k-1)(r_k-1)$.
 Again, since $\bask{k}$ is the widest $B$-path in $\Gamma$, we see similarly as before,
  \begin{equation*}
   \begin{split}
    |\ppot| & \leq  2 \sum_{i=1}^{l-1} r_i + r_l+ (\mu_t-1) r_k \leq \\
            & \leq  2 \sum_{i=1}^{k-1} r_i +(l-k-1)(r_k-1)+ (\mu_k-l+k+1) r_k \leq \\
            & \leq 2 \sum_{i=1}^{k-1} r_i + \mu_k r_k \; = \;  |\bask{k}| \, ,
   \end{split}
   \end{equation*}
 which proves the Lemma.
 \end{proof}
 
 \bigskip

\begin{theorem}\label{thm:premain}
 For $A \in \nb$ it follows that
  $$\nil{A} \leq \max\limits_{1 \leq i < t} 
    \{2i+\mu_{i+1}+\mu_{i+2}+\ldots+\mu_{i+r}; \; \mu_{i+1} - \mu_{i+r}\leq 1, \mu_i \ne \mu_{i+1}\} \, .$$
 \end{theorem}
 
 \medskip
 
 \begin{proof}
 Write $\sh{B}=(\mu_1,\mu_2,\ldots,\mu_t)$.
 We show by induction on $t$ that the longest possible path of $\nba$-digraphs is of the
 same length as a $B$-path.
 
 Suppose that the longest path of $\nba$-digraphs is included in digraph $\Gamma$ and
 denote it by $\ppot$. 
 Note that $\ppot$ contains vertices $(1,1)$ and $(1,\mu_1)$.
 If it does not, we can add an edge from $(1,1)$ to the first vertex of $\ppot$ or an edge
 from the last vertex of $\ppot$ to $(1,\mu_1)$ and lenghten it.
 
 \smallskip
 
 If $t=1$ then the longest path contains all vertices and therefore it is a $B$-path.
 
 Suppose that our claim holds for all partitions with at most $t-1$ parts.
 
 Fix an $x$, $1 \leq x \leq t$. If $\ppot$ does not contain any of the vertices $(x,y)$ for $y=1,2,\ldots,\mu_x$,
 then $\ppot$ is a path in a $({\mathcal N}_{B'},A)$-digraph, where 
 $$\sh{B'}=(\mu_1,\mu_2,\ldots,\mu_{x-1},\mu_{x+1},\mu_{x+2},\ldots,\mu_t).$$ By induction, there exists a $B'$-path 
 $\ppot_{B'}$ such that $|\ppot| \leq |\ppot_{B'}|$. If a $B'$-path is not already a $B$-path, 
 it can be lengthened to a $B$-path. Therefore there exists a $B$-path $\bas$ in $\Gamma$, such that 
 its length is greater than or equal to $|\ppot|$. By assumption that $\ppot$ is the longest path in $\Gamma$
 it follows that $|\ppot| = |\bas|$.
 
 Suppose now that for each $x$ the path $\ppot$ contains a vertex $(x,y_x)$ for some $y_x$.
 The basic idea of the proof is to show the following claim: 
 
 \emph{In $\Gamma$ there exists a path $\ppot'$ of the same length as $\ppot$ such that its first $t$ vertices are 
 $(1,y_1)$, $(2,y_2)$,..., $(t,y_t)$.} 
 
 Since $\ppot'$ is the longest path in $\Gamma$ it follows that $y_i=1$ for $i=1,2,\ldots,t$. 
 Otherwise the path that contains vertices $(1,1)$, $(2,1)$,..., $(t,1)$, $(t,2)$,..., $(t,y_t)$ and is after 
 $(t,y_t)$ equal to $\ppot'$ is longer than $\ppot'$.
 
 \smallskip
 
 Let us prove the claim. Let $y_t > 1$ be the smallest integer such that $\ppot$ contains vertex $(t,y_t)$. 
 Since $\ppot$ is the longest path, $(t,y_t)$ is not the last vertex of $\ppot$, and so for $x=1,\ldots,t$,
 there exists a vertex $(x,y^3_x)$ that $\ppot$ visits after $(t,y_t)$, i.e. $y^3_x \geq y_t$. 
 
 Next, suppose that there exists some $k$, where $1 \leq k < t$ and $y_k^1 < y_k^2 \leq y_t$, such that $\ppot$ 
 contains vertices $(k,y_k^1)$ and $(k,y_k^2)$ with the following properties:  there does 
 not exist a vertex $(k,y_k^0)$ of $\ppot$ such that $1\leq y_k^0 \leq y_k^1$ and $y_k^2$ is such that 
 there does not exist a vertex between $(k,y_k^2)$ and $(t,y_t)$ in $\ppot$ with $k$ as its first coordinate. 
 Denote by $\ppot_1$ (resp. $\ppot_2)$ the subpath of $\ppot$ with its initial vertex $(k,y_k^1)$ 
 (resp. $(k,y_k^2)$) and its final vertex $(k,y_k^2)$ (resp. $(k,y_k^3)$) and define 
 $\ppot_0=\ppot \, \backslash \, \{\ppot_1, \ppot_2\}$.
 In $\Gamma$ there exist the path $\ppot_1'$ parallel to 
 $\ppot_2$ starting at the $(k,y_k^1)$ and the path $\ppot_2'$ parallel to 
 $\ppot_1$ ending at the $(k,y_k^3)$. Then, the initial vertex of $\ppot'_2$ coincides with the final vertex of $\ppot'_1$.
 Thus $\ppot_0 \cup \ppot_1'\cup \ppot_2'$ is a path in $\Gamma$ of the same length as $\ppot$ and it contains only 
 one vertex with its first coordinate $k$ before $(t,y_t)$.
 
 By repeating this swap for $k=1,2,\ldots,t$, we obtain a path $\ppot'$ such that it does not contain vertices
 $(x,y_x^1)$ and $(x,y_x^2)$, with $y_x^1 < y_x^2 \leq y_t$ for all $x=1,2,\ldots, t-1$. 
 As argued above, this forces $\ppot'$ to be a path in the $\nba$-digraph $\Gamma$ that contains vertices 
 $(x,1)$ for all $x=1,2,\ldots,t$. By Lemma \ref{thm:path2} there exists a $B$-path $\bask{k}$ such that 
 $|\ppot|=|\ppot'| \leq 
 |\bask{k}|=2k+\mu_{k+1}+\mu_{k+2}+\ldots+\mu_{k+r}$ for some $r$, $\mu_k-\mu_{k+r}\leq 1$.
 \end{proof}

\bigskip

 To prove Theorem \ref{thm:main}, it only remains to show that the maximum of Theorem \ref{thm:premain} is attained. 

\bigskip
 
\begin{proof} (of {\bf Theorem \ref{thm:main}})
 Let $\sh{B}=(\mu_1,\mu_2,\ldots,\mu_t)=(m_1^{r_1},m_2^{r_2},\ldots,m_l^{r_l})$. We define the $B$-path $\ppot$, for which the maximum 
  $$\max\limits_{1 \leq i < t}\{2i+\mu_{i+1}+\mu_{i+2}+\ldots+\mu_{i+r}; \; \mu_{i+1} - \mu_{i+r}\leq 1; \; \mu_{i} \ne \mu_{i+1}\}$$
 is attained.
 
 Let $k$ and $s$ be such that 
 $\max\limits_{1 \leq i < t}\{2i+\mu_{i+1}+\mu_{i+2}+\ldots+\mu_{i+r}; \; \mu_{i+1} - \mu_{i+r}\leq 1; \; \mu_{i} \ne \mu_{i+1}\}=
 2k+\mu_{k+1}+\mu_{k+2}+\ldots+\mu_{k+s}$.
  
 \smallskip
 
 Denote by  $\ppot$ the $B_{(k+1)}$-path in $\Gamma$. 
 It follows that $|\ppot|=2k+\mu_{k+1}+\mu_{k+2}+\ldots+\mu_{k+s}$.  
 The path $\ppot$ is such that it can be completed (by drawing parallel edges) 
 to an $\nba$-digraph. (For examples of such $\nba$-digraphs see Example \ref{ex:2}.) 
 By the proof of Theorem \ref{thm:premain} it follows that $\ppot$ is a longest path in $\Gamma$.
 
 For a generic matrix $M_\Gamma$, corresponding to $\Gamma$, it follows that 
 $\nil{M}=|\ppot|$. 
 Take $A \in \nb$ to be a positive matrix with the same pattern 
 as $M$. 
 Since $\ch{\FF}=0$, it follows that $[A^h]_{ij} =0$ if and only if $[M^h]_{ij}=0$ for all $h$, $i$, $j$.
 Therefore $\nil{A}=\nil{M}=2k+\mu_{k+1}+\mu_{k+2}+\ldots+\mu_{k+s}$ which proves the theorem.
\end{proof}

\bigskip

\section{Examples} \label{sec:nilex}

\bigskip

\begin{example}\label{ex:2}
As in Examples \ref{ex:43221a} and \ref{ex:43221b}, we again consider the case $\sh{B}=(4,3,2^2,1)$. 
By Theorem \ref{thm:main} it 
follows that $\nil{\nb}=\max \{7,9,9\}=9$. 

Consider the following matrices:
$$
A_1=\left[\begin{array}{cccc|ccc|cc|cc|c}
0 & 0 & 0 & 0 & a & 0 & 0 & 0 & 0 & 0 & 0 & 0  \\
0 & 0 & 0 & 0 & 0 & a & 0 & 0 & 0 & 0 & 0 & 0  \\
0 & 0 & 0 & 0 & 0 & 0 & a & 0 & 0 & 0 & 0 & 0  \\
0 & 0 & 0 & 0 & 0 & 0 & 0 & 0 & 0 & 0 & 0 & 0  \\
\hline
0 & d & 0 & 0 & 0 & 0 & 0 & b & 0 & 0 & 0 & 0  \\
0 & 0 & d & 0 & 0 & 0 & 0 & 0 & b & 0 & 0 & 0  \\
0 & 0 & 0 & d & 0 & 0 & 0 & 0 & 0 & 0 & 0 & 0  \\
\hline
0 & 0 & 0 & 0 & 0 & 0 & 0 & 0 & 0 & c & 0 & 0  \\
0 & 0 & 0 & 0 & 0 & 0 & 0 & 0 & 0 & 0 & c & 0  \\
\hline
0 & 0 & 0 & 0 & 0 & e & 0 & 0 & 0 & 0 & 0 & 0  \\
0 & 0 & 0 & 0 & 0 & 0 & e & 0 & 0 & 0 & 0 & 0  \\
\hline
0 & 0 & 0 & 0 & 0 & 0 & 0 & 0 & 0 & 0 & 0 & 0  
\end{array}\right],$$

$$ A_2=\left[\begin{array}{cccc|ccc|cc|cc|c}
0 & f & 0 & 0 & a & 0 & 0 & 0 & 0 & 0 & 0 & 0  \\
0 & 0 & f & 0 & 0 & a & 0 & 0 & 0 & 0 & 0 & 0  \\
0 & 0 & 0 & f & 0 & 0 & a & 0 & 0 & 0 & 0 & 0  \\
0 & 0 & 0 & 0 & 0 & 0 & 0 & 0 & 0 & 0 & 0 & 0  \\
\hline
0 & d & 0 & 0 & 0 & 0 & 0 & b & 0 & 0 & 0 & 0  \\
0 & 0 & d & 0 & 0 & 0 & 0 & 0 & b & 0 & 0 & 0  \\
0 & 0 & 0 & d & 0 & 0 & 0 & 0 & 0 & 0 & 0 & 0  \\
\hline
0 & 0 & 0 & 0 & 0 & 0 & 0 & 0 & 0 & c & 0 & 0  \\
0 & 0 & 0 & 0 & 0 & 0 & 0 & 0 & 0 & 0 & c & 0  \\
\hline
0 & 0 & 0 & 0 & 0 & e & 0 & 0 & 0 & 0 & 0 & 0  \\
0 & 0 & 0 & 0 & 0 & 0 & e & 0 & 0 & 0 & 0 & 0  \\
\hline
0 & 0 & 0 & 0 & 0 & 0 & 0 & 0 & 0 & 0 & 0 & 0  
\end{array}\right] $$
and
$$
A_3= \left[\begin{array}{cccc|ccc|cc|cc|c}
0 & 0 & 0 & 0 & a & 0 & 0 & 0 & 0 & 0 & 0 & 0  \\
0 & 0 & 0 & 0 & 0 & a & 0 & 0 & 0 & 0 & 0 & 0  \\
0 & 0 & 0 & 0 & 0 & 0 & a & 0 & 0 & 0 & 0 & 0  \\
0 & 0 & 0 & 0 & 0 & 0 & 0 & 0 & 0 & 0 & 0 & 0  \\
\hline
0 & d & 0 & 0 & 0 & 0 & 0 & b & 0 & 0 & 0 & 0  \\
0 & 0 & d & 0 & 0 & 0 & 0 & 0 & b & 0 & 0 & 0  \\
0 & 0 & 0 & d & 0 & 0 & 0 & 0 & 0 & 0 & 0 & 0  \\
\hline
0 & 0 & 0 & 0 & 0 & 0 & 0 & 0 & 0 & c & 0 & 0  \\
0 & 0 & 0 & 0 & 0 & 0 & 0 & 0 & 0 & 0 & c & 0  \\
\hline
0 & 0 & 0 & 0 & 0 & e & 0 & 0 & 0 & 0 & 0 & f  \\
0 & 0 & 0 & 0 & 0 & 0 & e & 0 & 0 & 0 & 0 & 0  \\
\hline
0 & 0 & 0 & 0 & 0 & 0 & 0 & 0 & g & 0 & 0 & 0 
\end{array}\right],
$$

where $a,b,c,d,e,f,g \in \FF$ are nonzero.
The generic matrices $M_1$, $M_2$ and $M_3$ with the same patterns as $A_1$, $A_2$ and $A_3$
have the following digraphs
 
        \begin{figure}[htb]
         \begin{center}
          \includegraphics[height=8.2cm, width=12cm]{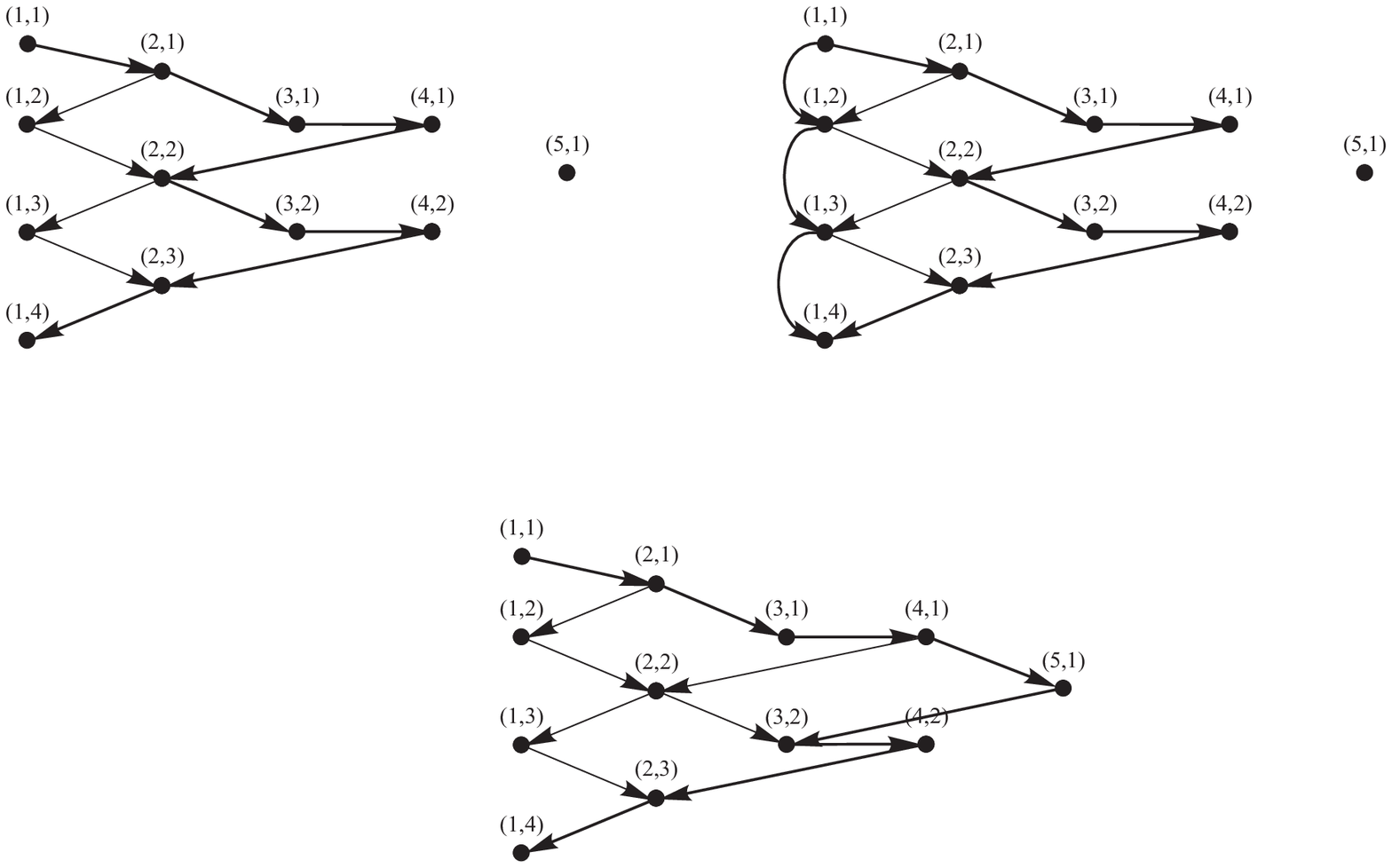}
         \end{center}
         \end{figure}

It can be easily seen that $\sh{A_1}=\sh{M_1}=(9,1^3)$, $\sh{A_2}=\sh{M_2}=(9,2,1)$
and $\sh{A_3}=\sh{M_3}=(9,3)$.
 \hfill$\square$
\end{example}

\bigskip

\begin{example}\label{ex:3}
Not all $\underline{\lambda}=(\lambda_1,\lambda_2,\ldots,\lambda_k) \in \partition{n}$, where $\lambda_1=\nil{\nb}$,
are in $\partition{\nb}$. For $\sh{B}=(6,4)$ it follows that $\nil{\nb}=6$. It can be shown that
$(6,3,1) \notin \partition{\nb}$.
  \hfill$\square$
\end{example}

 \bigskip

 \begin{example}
 Not all partitions in $\partition{\nb}$ can be obtained from a generic matrix corresponding to
 an $\nba$-digraph. 
 Let $\sh{B}=(5,3)$. It can be verified (for example, using {\tt Mathematica}) that 
 \begin{align*}
  \partition{\nb}=&
  \big\{(5,3), (5,2,1), (5,1^3), (4^2), (4,2^2), (4,2,1^2), (4,1^4), \\
     & \; (3^2,2), (3^2,1^2), (3,2^2,1), (3,2,1^3), (3,1^5), \\
     & \; (2^4),(2^3,1^2), (2^2,1^4), (2,1^6), (1^8)\big\}=\\
     = & \, \rp{5,3} \cup \big\{(4^2), (4,2^2), (4,2,1^2), (4,1^4), (3^2,1^2), (2^4) \big\}\, .
 \end{align*}
 One can check that there is no $\nba$-digraph $\Gamma$ with $\Delta(\Gamma)$ equal to 
 $(4^2)$, $(3^2,1^2)$ or $(2^4)$. However, there are matrices with these shapes in $\nb$.
 For example, for
 $$
 A_1=\left[\begin{array}{ccccc|ccc}
 0 & 1 & 0 & 0 & 0 &-1 & 0 & 0   \\
 0 & 0 & 1 & 0 & 0 & 0 &-1 & 0   \\
 0 & 0 & 0 & 1 & 0 & 0 & 0 &-1   \\
 0 & 0 & 0 & 0 & 1 & 0 & 0 & 0   \\
 0 & 0 & 0 & 0 & 0 & 0 & 0 & 0   \\
\hline
 0 & 0 & 1 & 0 & 0 & 0 &-1 & 0   \\
 0 & 0 & 0 & 1 & 0 & 0 & 0 &-1   \\
 0 & 0 & 0 & 0 & 1 & 0 & 0 & 0  
 \end{array}\right]$$
 and
 $$A_2=\left[\begin{array}{ccccc|ccc}
 0 & 1 & 0 & 0 & 0 & 1 & 0 & 0   \\
 0 & 0 & 1 & 0 & 0 & 0 & 1 & 0   \\
 0 & 0 & 0 & 1 & 0 & 0 & 0 & 1   \\
 0 & 0 & 0 & 0 & 1 & 0 & 0 & 0   \\
 0 & 0 & 0 & 0 & 0 & 0 & 0 & 0   \\
 \hline
 0 & 0 &-1 & 1 & 0 & 0 &-1 & 1   \\
 0 & 0 & 0 &-1 & 1 & 0 & 0 &-1   \\
 0 & 0 & 0 & 0 &-1 & 0 & 0 & 0  
 \end{array}\right]\, .
 $$
 
 we have $\sh{A_1}=(2^4)$, $\sh{A_2}=(3^2,1^2)$. 
 Let $M_i$, $i=1,2$, be generic matrices of the same shape as $A_i$. It can be proved, using
 Gansner-Saks Theorem, that $\sh{M_1}=\sh{M_2}=(5,3)$.
 \hfill$\square$
\end{example}


\bigskip
\bigskip




 
 \bigskip

\begin{thebibliography}{References}
  \bibitem{baranovsky}
   V.~Baranovsky: \emph{The variety of pairs of commuting nilpotent matrices is irreducible},
   Transform. Groups 6 (2001), no. 1, 3--8.
  \bibitem{basili}
   R.~Basili: \emph{On the irreducibility of commuting varieties of nilpotent matrices},
   J. Algebra 268 (2003), 58-80.
  \bibitem{brualdi}
   R.~A.~Brualdi, H.~J.~Ryser: \emph{Combinatorial matrix theory},
   Cambridge University Press, 1991.
  \bibitem{cmc}
   D.~H.~Collingwood, W.~M.~McGovern: \emph{Nilpotent orbits in semisimple Lie algebras},
   Van Nostrand Reinhold, 1993. 
  \bibitem{gansner}
   E.~R.~Gansner: \emph{Acyclic digraphs, Young tableaux and nilpotent matrices}, 
   SIAM J. Algebraic Discrete Methods 2 (1981), no. 4, 429--440.
  \bibitem{GLR}
   I.~Gohberg, P.~Lancaster, and L.~Rodman, \emph{Invariant Subspaces of Matrices with Applications}, 
   {Wiley-Interscience}, 1986.
  \bibitem{guralnick}
   R.~Guralnick: \emph{A note on commuting pairs of matrices}, Linear and Multilinear Algebra 31 (1992), no. 1-4, 71--75.
  \bibitem{gs}
   R.~Guralnick, B.~A.~Sethuraman: \emph{Commuting pairs and triples of matrices and related varieties}, 
   Linear Algebra Appl. 310 (2000), 139--148.
  \bibitem{gustafson}
   W.~H.~Gustafson: \emph{Modules and matrices}, Linear Algebra Appl. 157 (1991), 3-19.
  \bibitem{han}
   Y.~Han: \emph{Commuting triples of matrices}, Electron. J. Linear Algebra 13 (2005), 274-343. 
  \bibitem{ho}
   J.~Holbrook, M.~Omladi\v c: \emph{Approximating commuting operators}, Linear Algebra Appl. 327 (2001), no. 1-3, 131--149.
  \bibitem{kosiMST}
   T.~Ko\v sir: \emph{The Cayley-Hamilton theorem and inverse problems for multiparameter systems},
   Linear Algebra Appl. 367 (2003), 155-163.
  \bibitem{mt}
   T.~S.~Motzkin, O.~Taussky: \emph{Pairs of matrices with property ${\rm L}$}, Trans. Amer. Math. Soc. 73, (1952), 108--114.
  \bibitem{omladic}
   M.~Omladi\v c: \emph{A variety of commuting triples}, Linear Algebra Appl. 383 (2004), 233--245.
  \bibitem{saks}
   M.~Saks: \emph{Some sequences associated with combinatorial structures}, 
   Discrete Math. 59 (1986), no. 1-2, 135--166.  
  \bibitem{sivic}
   K.~\v Sivic: \emph{On varieties of commuting triples}, , Linear Algebra Appl. (2007), doi: 10.1016/j.laa.2007.11.004. 
\end{thebibliography}
\end{document}